\documentclass[reqno,12pt]{amsart}
\usepackage{amsmath,amsthm,amssymb,amsfonts,amscd,graphicx}
\input xy
\xyoption{all}
\theoremstyle{plain} 
	\newtheorem{thm}{Theorem}[section]
	\newtheorem*{thm*}{Theorem}
	\newtheorem{cor}[thm]{Corollary}
	\newtheorem{lem}[thm]{Lemma}
	
	\newtheorem{prop}[thm]{Proposition}
	
	\newtheorem*{conj*}{Conjecture}
	
\theoremstyle{definition}
	\newtheorem{defn}[thm]{Definition}

\theoremstyle{remark}
	\newtheorem{rem}[thm]{Remark}
	
	\newtheorem*{pf}{Proof}
\numberwithin{equation}{section}
\def\CC{{\mathbb C}}

\def\PP{{\mathbb P}}

\def\RR{{\mathbb R}}

\def\ZZ{{\mathbb Z}}

\def\D{{\mathcal D}}

\def\L{{\mathcal L}}

\def\N{{\mathcal N}}
\def\O{{\mathcal O}}

\def\S{{\mathcal S}}
\def\T{{\mathcal T}}

\newcommand{\per}{{\rm per}}

\def \mf#1#2#3#4{
\xymatrix{
{#1}\  \ar@<0.4ex>[r]^{{#2}} & \ {#4}
\ar@<0.4ex>[l]^{{#3}}
}
}

\def \mfs#1#2#3#4{\!
\xymatrix@C=1,5em{{#1} \! \ar@<0.2ex>[r]^{{#2}} & \! {#4}
\ar@<0.2ex>[l]^{{#3}}
}
\!}

\def \mfl#1#2#3#4{
\xymatrix@C=2.6em{{#1}\  \ar@<0.4ex>[r]^{{#2}} &\  {#4}
\ar@<0.2ex>[l]^{{#3}}
}
}

\def \mfss#1#2#3#4{\!
\xymatrix@C=1.5em{{#1} \ar@<0.3ex>[r]^{{#2}} & {#4}
\ar@<0.3ex>[l]^{{#3}}
}
\!}
\begin{document}
\title{A note on entropy of auto-equivalences:
lower bound and the case of orbifold projective lines}
\date{\today}
\author{Kohei Kikuta}
\address{Department of Mathematics, Graduate School of Science, Osaka University, 
Toyonaka Osaka, 560-0043, Japan}
\email{k-kikuta@cr.math.sci.osaka-u.ac.jp}
\author{Yuuki Shiraishi}
\address{Department of Mathematics, Graduate School of Science, Kyoto University, 
Kyoto 606-8502, Japan}
\email{yshiraishi@math.kyoto-u.ac.jp}
\author{Atsushi Takahashi}
\address{Department of Mathematics, Graduate School of Science, Osaka University, 
Toyonaka Osaka, 560-0043, Japan}
\email{takahashi@math.sci.osaka-u.ac.jp}
\begin{abstract}
Entropy of categorical dynamics is defined by Dmitrov--Haiden--Katzarkov--Kontsevich. 
Motivated by the fundamental theorem of the topological entropy due to Gromov--Yomdin, 
it is natural to ask an equality between the entropy and the spectral radius of induced morphisms on the numerical Grothendieck group. 

In this paper, we add two results on this equality: the lower bound in a general setting and the equality for orbifold projective lines. 
\end{abstract}
\maketitle
\markboth{KOHEI KIKUTA, YUUKI SHIRAISHI, AND ATSUSHI TAKAHASHI}{A NOTE ON ENTROPY OF AUTO-EQUIVALENCES}
\section{Introduction}
It is interesting to bring some dynamical view points into the category theory. 
Motivated by the classical theory of dynamical systems, 
the notion of entropy of categorical dynamical systems (entropy of endo-functors for short) is 
defined by Dimitrov--Haiden--Katzarkov--Kontsevich \cite{DHKK}. 
The entropy of endo-functors is actually similar to the topological entropy in the sense of sharing many properties (Lemma~\ref{entropy}, \ref{admissible}, \ref{ess-surj}). 
Moreover, the entropy of the derived pull-back of a surjective endomorphism of a smooth projective variety over $\CC$ is equal to its topological entropy \cite{KT}.
In other words, the entropy of endo-functors can be thought of as a categorical generalization of the topological entropy. 

In this paper, we add two results on the entropy of endo-functors.
The first one is that, for the perfect derived categories $\per(B)$ of a smooth proper differential graded algebra $B$, the lower bound of the entropy $h(F)$ of an endo-functor $F$ is given by the natural logarithm of the spectral radius $\rho(\N(F))$ on the numerical Grothendieck group, 
called the (numerical) Gromov--Yomdin type inequality (See also \cite[Conjecture~5.3]{KT}).
It is motivated by the fundamental theorem of the topological entropy for complex dynamics on algebraic varieties due to Gromov--Yomdin \cite{Gro1,Gro2,Yom}:
\begin{thm*}[Theorem~\ref{nGro-lower}]
For each endo-functor $F$ of $\per(B)$ admitting left or right adjoint functors, such that $F^nB\not\cong 0$ for $n\ge 0$, 
we have 
\begin{equation}\label{eq-lower}
h(F)\geq\log\rho(\N(F)).
\end{equation}
\end{thm*} 
For the proof, we use some norm inspired by the theory of dynamical degree and algebraic cycles due to Truong~\cite{Tru}. 
Ikeda shows this inequality by the mass-growth for Bridgeland's stability conditions \cite{Ike}.  

The equality in the Gromov--Yomdin type inequality is now known to hold for elliptic curves \cite{Kik}, varieties with the ample (anti-)canonical sheaf \cite{KT} and abelian surfaces \cite{Yos}, which gives some applications to the topological entropy of dynamics on moduli spaces of stable objects in the sense of Bridgeland \cite{Ouc1,Yos}. 
But, in general, it does not hold for some Calabi-Yau varieties \cite{Fan,Ouc2}. 
As a corollary of the first main theorem, it is easy to show the equality for derived categories of hereditary finite dimensional algebras (Proposition~\ref{GY-hereditary}, Corollary~\ref{GY-Dynkin}).

The second result of  this paper claims the equality for the derived category $\D^b(\PP^1_{A,\Lambda})$ of an orbifold projective line $\PP^1_{A,\Lambda}$ introduced by Geigle--Lenzing \cite{GL}. 
Orbifold projective lines are important and interesting objects since they are not only 
in the next class to hereditary finite dimensional algebras but few examples 
whose homological and classical mirror symmetry are well-understood (cf. \cite{IST1,IST2,IT,Kea,Ros,ST,Tak1,Tak2,Ued}): 
\begin{thm*}[Theorem~\ref{main2}]
For each auto-equivalence $F$ of $\D^b(\PP^{1}_{A,\Lambda})$, we have
\begin{equation}
h(F)=\log\rho(\N(F)). 
\end{equation}
Moreover, $\rho(\N(F))$ is an algebraic number and $h(F)=0$ if $\chi_{A}\ne 0$.
\end{thm*}
It is an important and interesting problem to find a characterization of endo-functors attaining the lower bound of the inequality (\ref{eq-lower}). 

\begin{sloppypar}
{\bf Acknowledgements}.\   
The first named author is supported by JSPS KAKENHI Grant Number JP17J00227. 
The second named author is supported by Research Fellowship of Japan Society for the Promotion for Young Scientists. 
The third named author is supported by JSPS KAKENHI Grant Number JP16H06337, JP26610008. 
\end{sloppypar}
\section{Preliminaries}
\subsection{Notations and terminologies}
Throughout this paper, we work over the base field $\CC$ and all triangulated categories are $\CC$-linear and not equivalent to the zero category. 
The translation functor on a triangulated category is denoted by $[1]$.
All (triangulated) functors are $\CC$-linear. 

A triangulated category $\T$ is called {\it split-closed} if every idempotent in $\T$ splits, namely, if it contains all direct summands of its objects, and it is called {\it thick} if it is split-closed and closed under isomorphisms.
For an object $M\in\T$, we denote $\langle M\rangle$ by the smallest thick triangulated subcategory containing $M$.
An object $G\in\T$ is called a {\it split-generator} if $\langle G\rangle=\T$.
A triangulated category $\T$ is said to be {\it of finite type} if for all $M,N\in\T$ we have $\sum_{n\in\ZZ}\dim_\CC{\rm Hom}_{\T}(M,N[n])<\infty$.

\subsection{Complexity}
From now on, $\T,\T'$ denote triangulated categories of finite type.
\begin{defn}[Definition~2.1 in \cite{DHKK}]
For each $M,N\in \T$, define the function $\delta_{\T,t}(M,N):\RR\longrightarrow\RR_{\geq0}\cup\{ \infty \}$ in $t$ by
{\small 
\begin{equation*}
\delta_{\T,t}(M,N):= 
\begin{cases}
0
 & \text{ if }N\cong0\\
\inf\left\{
\displaystyle\sum_{i=1}^p {\rm exp}(n_i t)~
\middle
|~
\begin{xy}
(0,5) *{0}="0", (20,5)*{A_{1}}="1", (30,5)*{\dots}, (40,5)*{A_{p-1}}="k-1", (60,5)*{N\oplus N'}="k",
(10,-5)*{M[n_{1}]}="n1", (30,-5)*{\dots}, (50,-5)*{M[n_{p}]}="nk",
\ar "0"; "1"
\ar "1"; "n1"
\ar@{.>} "n1";"0"
\ar "k-1"; "k" 
\ar "k"; "nk"
\ar@{.>} "nk";"k-1"
\end{xy}
\right\}
 & \text{ if }N\in\langle M\rangle \\
\infty
 & \text{ if }N\not\in\langle M\rangle.
\end{cases}
\end{equation*}
}
The function $\delta_{\T,t}(M,N)$ is called the {\it complexity} of $N$ with  respect to $M$.
\end{defn}
\begin{rem}
If $\T$ has a split-generator $G$ and $M\in\T$ is not isomorphic to a zero object, 
then an inequality $1\leq\delta_{\T,0}(G,M)<\infty$ holds.
\end{rem} 
We recall some basic properties of the complexity.  
\begin{lem}
Let $M_1,M_2,M_3,M_4\in\T$.
\begin{enumerate}
\item
If $M_2\in\langle M_1\rangle$ and $M_2\not\cong0$, then $0<\delta_{\T,t}(M_1,M_2)$. 
\item
If $M_1\cong M_3$, then $\delta_{\T,t}(M_1,M_2)=\delta_{\T,t}(M_3,M_2)$. 
\item
If $M_2\cong M_3$, then $\delta_{\T,t}(M_1,M_2)=\delta_{\T,t}(M_1,M_3)$. 
\item
If $M_2\not\cong0$, then $\delta_{\T,t}(M_1,M_3)\leq\delta_{\T,t}(M_1,M_2)\delta_{\T,t}(M_2,M_3)$.
\item
We have $\delta_{\T,t}(M_4,M_2)\leq\delta_{\T,t}(M_4,M_1)+\delta_{\T,t}(M_4,M_3)$ for an exact triangle $M_1\to M_2\to M_3$.
\item
We have $\delta_{\T',t}(F(M_1),F(M_2))\leq\delta_{\T,t}(M_1,M_2)$ for any triangulated functor $F:\T\longrightarrow \T'$.
\end{enumerate}
\end{lem}
\begin{lem}
Let $\D^b(\CC)$ be the bounded derived category of finite dimensional $\CC$-vector spaces. 
For $M\in\D^b(\CC)$, we have the following inequality 
\begin{equation}\label{eq:2}
\delta_{\D^b(\CC),t}(\CC,M)=\sum_{l\in\ZZ} \left(\dim_{\CC}H^{l}(M)\right)\cdot e^{-lt}.
\end{equation}
\end{lem}

\subsection{Entropy of endo-functors}\label{subsection}
Endo-functor $F$ means triangulated functor $F:\T\to\T$. 
We assume that all endo-functors of $\T$ satisfy that $F^n G\not\cong0$ for $n\ge0$ (if $\T$ has a split-generator $G$). 
\begin{defn}[Definition~2.4 in \cite{DHKK}]
Let $G$ be a split-generator of $\T$ and $F$ an endo-functor of $\T$. 
The {\it entropy} of $F$ is the function $h_t(F):\RR\longrightarrow\{ -\infty \}\cup\RR$ given by
\begin{equation}
h_{t}(F):=\displaystyle\lim_{n\rightarrow\infty}\frac{1}{n}\log \delta_{\T,t}(G,F^{n}G).
\end{equation}
\end{defn}
It follows from \cite[Lemma~2.5]{DHKK} that the entropy is well-defined and doesn't depend on the choice of split-generators.
\begin{lem}\label{two-split-gen}
Let $G,G'$ be split-generators of $\T$ and $F$ an endo-functor of $\T$. 
The entropy $h_t(F)$ of $F$ is given by
\begin{equation}
h_t(F)=\displaystyle\lim_{n\rightarrow\infty}\frac{1}{n}\log \delta_{\T,t}(G,F^{n}G').
\end{equation}
\end{lem}
The three lemmas below show that the entropy of endo-functors is similar to the topological entropy. 
\begin{lem}\label{entropy}
Let $G$ be a split-generator of $\T$ and $F_1,F_2$ endo-functors of $\T$. 
\begin{enumerate}
\item
If $F_{1} \cong F_{2}$, then $h_t(F_{1})=h_t(F_{2})$. 
\item
We have $h_t(F_{1}^m)= mh_t(F_{1})$ for $m\geq1$.
\item
we have $h_t(F_1F_2)=h_t(F_2F_1)$.
\item
If $F_1F_2\cong F_2F_1$, then $h_t(F_1F_2)\leq h_t(F_1)+h_t(F_2)$.
\item
If $F_1=F_2[m]~(m\in\ZZ)$, then $h_t(F_1)=h_t(F_2)+mt$. 
\end{enumerate}
\end{lem}
\begin{lem}\label{admissible}
Let $F_i$ be an endo-functor of $\T_i$ with a split-generator $G_i~(i=1,2)$.
If there exists a fully faithful functor $F':\T_2\to\T_1$, which has left and right adjoint functors, such that $F'F_2\simeq F_1F'$, then $h_t(F_2)\leq h_t(F_1)$. 
\end{lem}
\begin{lem}\label{ess-surj}
Let $F_i$ be an endo-functor of $\T_i$ with a split-generator $G_i~(i=1,2)$.
If there exists a essentially surjective functor $F':\T_1\to\T_2$ such that $F'F_1\simeq F_2F'$, then $h_t(F_2)\leq h_t(F_1)$.  
\end{lem}
As a corollary of Lemma~\ref{ess-surj}, we have the following 
\begin{cor}
Let $F'$ an auto-equivalence of $\T$. 
The entropy is a class function, namely, $h_t(F'FF'^{-1})=h_t(F)$. 
\end{cor}

Let $B$ be a smooth proper differential graded (dg) $\CC$-algebra $B$ and  
$\per(B)$ the perfect derived category of dg $B$-modules, 
the full triangulated subcategory of the derived category $\D(B)$ of dg $B$-modules containing $B$
closed under isomorphisms and taking direct summands. 
By definition, $B$ is a split-generator of $\per(B)$.

The following proposition enables us to compute entropy. 
\begin{prop}[Theorem~2.7 in \cite{DHKK}]\label{smooth-proper}
Let $G,G'$ be split-generators of $\per(B)$ and $F$ an endo-functor of $\per(B)$. 
The entropy $h_t(F)$ is given by
\begin{equation}
h_t(F)=\lim_{n\rightarrow\infty}\frac{1}{n}\log\delta'_{\per(B),t}(G,F^nG'),
\end{equation}
where 
\begin{equation}
\delta'_{\per(B),t}(M,N):=\sum_{m\in\ZZ}\dim_{\CC} {\rm Hom}_{\per(B)}(M,N[m]) e^{-mt},\quad M,N\in\per(B).
\end{equation}
\end{prop}
\begin{pf}
The following is proven in the proof of \cite[Theorem~2.7]{DHKK}.
\begin{lem}\label{delta-delta'}
For each $M\in\per(B)$, there exist $C_1(t),C_2(t)$ for $t\in\RR$ such that 
\[
C_1(t)\delta_{\per(B),t}(G,M)
\leq\delta'_{\per(B),t}(G,M)
\leq C_2(t)\delta_{\per(B),t}(G,M).
\]
In particular, for each $M\in\per(B)$ we have
\begin{equation}
\lim_{n\rightarrow\infty}\frac{1}{n}\log\delta_{\per(B),t}(G,M)
=\lim_{n\rightarrow\infty}\frac{1}{n}\log\delta'_{\per(B),t}(G,M).
\end{equation}
\end{lem}
Together with Lemma~\ref{two-split-gen}, we have
\[
h_t(F)=\lim_{n\rightarrow\infty}\frac{1}{n}\log\delta_{\per(B),t}(G,F^nG')
=\lim_{n\rightarrow\infty}\frac{1}{n}\log\delta'_{\per(B),t}(G,F^nG').
\]
We finished the proof of the proposition.
\qed
\end{pf}
In order to state the first main theorem, we prepare some terminologies.
For $M,N\in\per(B)$, set 
\begin{equation}
\chi(M,N):=\sum_{n\in\ZZ}(-1)^n\dim_\CC{\rm Hom}_{\per(B)}(M,N[n]).
\end{equation}
It naturally induces a bilinear form on the Grothendieck group $K_0(\per(B))$ of $\per(B)$, called the {\it Euler form}, 
which is denoted by the same letter $\chi$. 
Then the {\it numerical Grothendieck group} $\N(\per(B))$ is defined as the quotient of $K_0(\per(B))$ by the radical of $\chi$ 
(which is well-defined by the Serre duality). 
It is important to note that $\N(\per(B))$ is a free abelian group of finite rank by Hirzebruch-Riemann-Roch theorem \cite{Shk,Lun}.
If an endo-functor $F$ of $\per(B)$ admits left or right (hence both by the Serre duality) adjoint functors, it respects the radical of $\chi$. 
Therefore, it induces an endomorphism $\N(F)$ on $\N(\per(B))$. 
Note that an endo-functor lifting to a dg endo-functor of the dg category $\per_{dg}(B)$ admits adjoint functors. 
The {\it spectral radius} $\rho(\N(F))$ of $\N(F)$ is the maximum of absolute values of eigenvalues of $\CC$-linear endomorphism $\N(F)\otimes_\ZZ \CC$.
Set $\delta_\T:=\delta_{\T,0},\delta'_\T:=\delta'_{\T,0},h:=h_0$.

Inspired by the theory of  dynamical degree and algebraic cycles due to Truong~(cf. \cite[eq. (3.2)]{Tru}), we show the following:
\begin{thm}\label{nGro-lower}
For each endo-functor $F$ of $\per(B)$ admitting left or right adjoint functors, we have 
\begin{equation}
h(F)\geq\log\rho(\N(F)).
\end{equation}
\end{thm}
\begin{pf}
Let $v_1,\cdots,v_p~(v_i=[M_i],M_i\in\per(B))$ be a fixed basis of $\N(\per(B))$.  
Set $M_0:=\oplus_i M_i,\N(\per(B))_{\RR}:=\N(\per(B))\otimes_{\ZZ}\RR,\N(F)_{\RR}:=\N(F)\otimes_{\ZZ}\RR,\chi_{\RR}:=\chi\otimes_\ZZ\RR$.
Define a norm $\|\cdot\|$ on $\N(\per(B))_\RR$ by 
\begin{equation}
\|v\|:=\displaystyle\sum_{i=1}^p|\chi_{\RR}(v_i,v)|,\quad v\in\N(\per(B))_\RR,
\end{equation}
which induces an operator norm of $\N(F)_{\RR}$, that is, $\|\N(F)_{\RR}\|:=\displaystyle\sup_{\|v\|=1}\|\N(F)_{\RR}v\|$. 
By the compactness of the subset $\{\|v\|=1\}\subset\N(\per(B))_{\RR}$, 
there exists a positive number $C>0$ such that 
\[
\displaystyle\sum_{i,j=1}^p|\chi(v_i,\N(F)v_j)|\geq C\cdot\|\N(F)_{\RR}\|.
\]
Note that $B\oplus M_0$ is a split-generator of $\per(B)$. 
By Proposition \ref{smooth-proper}, the statement follows from
\begin{eqnarray*}
h(F)
&=&\displaystyle\lim_{n\to\infty}\frac{1}{n}\log\delta'_{\per(B)}(B\oplus M_0,F^n(B\oplus M_0))\\
&\geq&\displaystyle\lim_{n\to\infty}\frac{1}{n}\log\delta'_{\per(B)}(M_0,F^nM_0)\\
&\geq&\displaystyle\lim_{n\to\infty}\frac{1}{n}\log\sum_{i,j}|\chi(v_i,\N(F^n)v_j)|\\
&\geq&\displaystyle\lim_{n\to\infty}\frac{1}{n}\log\|\N(F)^n_{\RR}\|=\log\rho(\N(F)).
\end{eqnarray*}
\qed
\end{pf}
Let ${\rm Auteq}(\T)$ be the group of (natural isomorphism classes of) auto-equivalences of a triangulated category $\T$. 
\begin{prop}\label{GY-hereditary}
Let $B$ be a hereditary finite dimensional $\CC$-algebra.
For each auto-equivalence $F\in{\rm Auteq}(\per(B))$, we have
\begin{equation}
h(F)=\log\rho(\N(F)).
\end{equation}
\end{prop}
\begin{pf}
Due to Theorem~\ref{nGro-lower}, we only need to show the upper bound.
Let $P_1,\dots,P_{\dim_{\CC} B}$ be indecomposable modules. 
Each auto-equivalence $F$ sends an indecomposable object to an indecomposable one. 
Since $B$ is hereditary, there exists $m\in \ZZ$ such that the indecomposable object $F^{n}(P_{i})[m]$ 
is isomorphic to an object concentrated in degree zero, namely, a $B$-module.
By Proposition~\ref{smooth-proper}, we have 
\begin{eqnarray*}
h(F)
&=& \lim_{n\rightarrow\infty}\frac{1}{n}\log\delta'_{\per(B)}(B,F^{n}(\oplus_{i=1}^{\dim_{\CC}B}P_{i}))\\
&=& \lim_{n\rightarrow\infty}\frac{1}{n}\log\sum_{i=1}^{\dim_\CC B}\left|\chi(B,F^{n}(P_{i}))\right| \le\log \rho(\N(F)).
\end{eqnarray*}
\qed
\end{pf}
\begin{cor}\label{GY-Dynkin}
Suppose that $B=\CC\vec{\Delta}$ for some Dynkin quiver $\vec{\Delta}$. 
Then, we have 
\begin{equation}
h(F)=\log\rho(\N(F))=0.
\end{equation}
\end{cor}
\begin{pf}
It is known by \cite[Theorem~3.8]{MY}, that
\begin{equation}
{\rm Auteq}(\per(B))\cong \langle \S_{B}, \S_{B}[-1] \rangle \times {\rm Aut}(\vec{\Delta}),
\end{equation}
where $\S_B$ is the Serre functor of $\per(B)$
and ${\rm Aut}(\vec{\Delta})$ is the finite subgroup of ${\rm Auteq}(\per(B))$
consisting of automorphisms of $\vec{\Delta}$.
Again, by \cite[Theorem~3.8]{MY}, $\S_{B}$ is of finite order up to translation. 
The statement follows from Lemma~\ref{entropy} {\rm (ii)}, {\rm (iv)} and {\rm (v)}. 
\qed
\end{pf}
\section{Orbifold projective lines}
In this section, we shall show the Gromov--Yomdin type theorem for the entropy of an auto-equivalence on
the derived category $\D^{b}(\PP^{1}_{A, \Lambda})$ of coherent sheaves on an orbifold projective line $\PP^{1}_{A, \Lambda}$.   
We first recall the definition of orbifold projective line in \cite{GL}.

Let $r\ge 3$ be a positive integer.
Let $A = (a_{1},\ldots,a_{r})$ be a multiplet of positive integers and 
$\Lambda=(\lambda_1,\ldots,\lambda_r)$ a multiplet of pairwise distinct elements of $\PP^{1}(\CC)$ normalized such that 
$\lambda_1=\infty$, $\lambda_2=0$ and $\lambda_3=1$. 
In order to introduce an orbifold projective line,
we prepare some notations.
\begin{defn}\label{orb proj line}
Let $r$, $A$ and $\Lambda$ be as above.
\begin{enumerate}
\item
Define a ring $R_{A,\Lambda}$ by 
\begin{subequations}
\begin{equation} 
R_{A,\Lambda}:=\CC[X_1,\dots,X_r]\left/I_{\Lambda}\right.,
\end{equation}
where  $I_\Lambda$ is an ideal generated by $r-2$ homogeneous polynomials
\begin{equation}
X_i^{a_i}-X_2^{a_2}+\lambda_{i} X_1^{a_1},\quad i=3,\dots, r.
\end{equation}
\end{subequations}
\item
Denote by $L_A$ an abelian group generated by $r$-letters $\vec{X_i}$, $i=1,\dots ,r$ defined as the quotient 
\begin{subequations}
\begin{equation}
L_A:=\bigoplus_{i=1}^r\ZZ\vec{X}_{i}\left/M_A\right. ,
\end{equation}
where  $M_A$ is the subgroup generated by the elements
\begin{equation}
a_i\vec{X}_i-a_j\vec{X}_j,\quad 1\le i<j\le r.
\end{equation}
\end{subequations}
\item
Set 
\begin{equation}
a:={\rm l.c.m}(a_{1},\dots,a_{r}),\quad
\mu_{A}:=2+\sum^{r}_{i=1}(a_{i}-1), \quad \chi_{A}:=2+\sum^{r}_{i=1}(\frac{1}{a_{i}}-1).
\end{equation}
\end{enumerate}
\end{defn}
We then consider the following quotient stack$:$
\begin{defn}\label{defn:gl}
Let $r$, $A$ and $\Lambda$ be as above.
Define a stack $\PP^1_{A,\Lambda}$ by
\begin{equation}
\PP^1_{A,\Lambda}:=\left[\left({\rm Spec}(R_{A,\Lambda})\backslash\{0\}\right)/{\rm Spec}({\CC L_A})\right],
\end{equation}
which is called the {\it orbifold projective line} of type $(A,\Lambda)$. 
\end{defn}
The orbifold projective line is a Deligne--Mumford stack whose coarse moduli space is 
a smooth projective line $\PP^1$.

Denote by ${\rm gr}^{L_{A}}(R_{A,\Lambda})$ the abelian category of finitely generated $L_{A}$-graded $R_{A,\Lambda}$-modules 
and denote by ${\rm tor}^{L_{A}}(R_{A,\Lambda})$ the full subcategory of ${\rm gr}^{L_{A}}(R_{A,\Lambda})$ 
whose objects are finite-dimensional $L_{A}$-graded $R_{A,\Lambda}$-modules. 
It is known (cf. \cite[Section~1.8]{GL}) that the abelian category ${\rm Coh}(\PP^1_{A,\Lambda})$ of coherent sheaves is given by
\begin{equation}
{\rm Coh}(\PP^1_{A,\Lambda})={\rm gr}^{L_{A}}(R_{A,\Lambda})/{\rm tor}^{L_{A}}(R_{A,\Lambda}).
\end{equation}
Denote by $\D^b(\PP^1_{A,\Lambda})$ the bounded derived category $\D^b({\rm Coh}(\PP^1_{A,\Lambda}))$ of ${\rm Coh}(\PP^1_{A,\Lambda})$.

For each $\vec{l}\in L_{A}$, set
\begin{equation}
\O_{\PP^{1}_{A,\Lambda}}(\vec{l}):=[R_{A,\Lambda}(\vec{l})]\in {\rm Coh}(\PP^1_{A,\Lambda}),
\end{equation}
where $(R_{A,\Lambda}(\vec{l}))_{\vec{l'}}:=(R_{A,\Lambda})_{\vec{l}+\vec{l'}}$. 

Set $\vec{c}:=a_{1}\vec{x}_{1}=\dots=a_{r}\vec{x}_{r}$. 
The element $\vec{x}\in L_{A}$ has the unique expression of the form
\begin{equation}
\vec{x}=l\vec{c}+\sum_{i=1}^{r}p_{i}\vec{x}_{i}, \quad \ 0\le p_{i}\le a_{i}-1. 
\end{equation}
We say that $\vec{x}$ is {\it positive} if $\vec{x}\ne 0$, $l\ge 0$ and $x_{i}\ge 0$ for $i=1,\dots,r$.

For a $\CC$-module $M$, set $M^{*}:={\rm Hom}_{\CC}(M, \CC)$.
\begin{prop}[Section 1.8.1 and Section 2.2 in \cite{GL}]\label{hom-str-of-line}
We have the following:
\begin{enumerate}
\item For $\vec{x}, \vec{y}\in L_{A}$ with $\vec{x}-\vec{y}$ positive,
\begin{equation}
{\rm Hom}(\O_{\PP^1_{A,\Lambda}}(\vec{x}),\O_{\PP^1_{A,\Lambda}}(\vec{y}))=0.  
\end{equation} 
\item Set $\vec{\omega}:=(r-2)\vec{c}-\sum_{i=1}^{r}\vec{x_{i}}\in L_{A}$.
For $M_{1}, M_{2}\in {\rm Coh}(\PP^1_{A,\Lambda})$, we have the Serre duality isomorphism: 
\begin{equation}
{\rm Ext}^{1}(M_{2}, M_{1})\cong {\rm Hom}(M_{1}, M_{2}\otimes \O_{\PP^{1}_{A,\Lambda}}(\vec{\omega}))^{*}.
\end{equation}
\item The category ${\rm Coh(\PP^1_{A,\Lambda})}$ is hereditary, namely, ${\rm Ext}^{i}(M_{1}, M_{2})=0$ for 
$M_{1}, M_{2}\in {\rm Coh}(\PP^1_{A,\Lambda})$ if $i\ne 0,1$.
\end{enumerate}
\end{prop}

\begin{rem}
It follows from Proposition~\ref{hom-str-of-line} {\rm (iii)} that
each indecomposable object of $\D^{b}(\PP^{1}_{A,\Lambda})$ is of the form $M[n]$ for 
some $M\in {\rm Coh}(\PP^{1}_{A,\Lambda})$ and $n\in \ZZ$.
\end{rem}

\begin{prop}[Section 1.8.1 and Section 4.1 in \cite{GL}]\label{generators-orb-proj-line}
The following sequences are full strongly exceptional collections:
\begin{eqnarray*}
(E_{1},\dots,E_{\mu_{A}}):=(\O_{\PP^1_{A,\Lambda}}, \O_{\PP^1_{A,\Lambda}}(\vec{x}_{1}),\dots,
\O_{\PP^1_{A,\Lambda}}((a_{1}-1)\vec{x}_{1}),\dots,\\
\O_{\PP^1_{A,\Lambda}}(\vec{x}_{r}),\dots,
\O_{\PP^1_{A,\Lambda}}((a_{r}-1)\vec{x}_{r}), \O_{\PP^{1}_{A,\Lambda}}(\vec{c}))
\end{eqnarray*}
\begin{eqnarray*}
(E^{*}_{\mu_{A}}, \dots E^{*}_{1}):=(
\O_{\PP^{1}_{A,\Lambda}}(-\vec{c}),\O_{\PP^1_{A,\Lambda}}(-(a_{r}-1)\vec{x}_{r}),\dots, 
\O_{\PP^1_{A,\Lambda}}(-\vec{x}_{r}),\dots,\\
\O_{\PP^1_{A,\Lambda}}(-(a_{1}-1)\vec{x}_{1}), \dots, \O_{\PP^1_{A,\Lambda}}(-\vec{x}_{1}), \O_{\PP^1_{A,\Lambda}}).
\end{eqnarray*}
In particular, $G:=\displaystyle\bigoplus_{i=1}^{\mu_{A}}E_{i}$ and $G^{*}:=\displaystyle\bigoplus_{i=1}^{\mu_{A}}E^{*}_{i}$ 
are split-generators of $\D^b(\PP^{1}_{A,\Lambda})$.
\end{prop}

It follows from Proposition~\ref{generators-orb-proj-line} that $\D^b(\PP^1_{A,\Lambda})\cong \per({\rm End}(G))$
since the triangulated category $\D^b(\PP^1_{A,\Lambda})$ is algebraic.
Denote by $\N(\PP^1_{A,\Lambda})$ its numerical Grothendieck group.

\begin{defn}[Section 2.5 in \cite{GL}]
Take $[1: \lambda ]\in \PP^{1}\setminus \{\lambda_{1},\dots,\lambda_{r}\}$. 
Define $S$ and $S_{i,j}$ for $i=1,\dots, r$ and $j=0, \dots, a_{i}-1$ by the following exact sequences:
\begin{equation}
0\rightarrow \O_{\PP^{1}_{A,\Lambda}}(\vec{0})\xrightarrow{X^{a_{1}}_{1}-\lambda X^{a_{2}}_{2}} \O_{\PP^{1}_{A,\Lambda}}(\vec{c}) 
\rightarrow S \rightarrow 0,
\end{equation}
\begin{equation}
0\rightarrow \O_{\PP^{1}_{A,\Lambda}}(j\vec{x}_{i})\xrightarrow{X_{i}} \O_{\PP^{1}_{A,\Lambda}}((j+1)\vec{x}_{i}) 
\rightarrow S_{i,j} \rightarrow 0.
\end{equation}
\end{defn}

\begin{defn}[Section 1.8.2 and Section 2.8 in \cite{GL}]
The {\it rank} and {\it degree} are homomorphisms defined as follows: 
\begin{equation}
{\rm rk}: \N(\PP^{1}_{A,\Lambda})\rightarrow \ZZ, \quad {\rm rk}([\O_{\PP^{1}_{A,\Lambda}}]):=1, \ \ 
{\rm rk}([S]):=0 \ \ \text{and} \ \ {\rm rk}([S_{i,j}]):=0,
\end{equation}
\begin{equation}
{\rm deg}: \N(\PP^{1}_{A,\Lambda})\rightarrow \ZZ, \quad {\rm deg}([\O_{\PP^{1}_{A,\Lambda}}(\vec{x_{i}})]):=\frac{a}{a_{i}}
\ \ \text{and} \ \ {\rm deg}([\O_{\PP^{1}_{A,\Lambda}}]):=0. 
\end{equation}
\end{defn}
\begin{defn}
Denote by ${\rm Pic}(\PP^{1}_{A,\Lambda})$ the group consisting of (isomorphism classes of) indecomposable objects
in ${\rm Coh}(\PP^{1}_{A,\Lambda})$ of rank one with multiplication
induced by the tensor product.
\end{defn}

\begin{lem}[Section 2.1 in \cite{GL}]
There is an isomorphism of abelian groups
\begin{equation}
L_{A}\cong {\rm Pic}(\PP^{1}_{A,\Lambda}), \quad \vec{x}_{i}:=[\vec{X}_{i}]\mapsto \O_{\PP^1_{A,\Lambda}}(\vec{x}_{i}).
\end{equation}
\end{lem}

One of our results is the following Gromov--Yomdin type theorem for an orbifold projective line:
\begin{thm}\label{main2}
For each auto-equivalence $F$ of $\D^b(\PP^{1}_{A,\Lambda})$, we have
\begin{equation}
h(F)=\log\rho(\N(F)). 
\end{equation}
Moreover, $\rho(\N(F))$ is an algebraic number and $h(F)=0$ if $\chi_{A}\ne 0$.
\end{thm}

\subsection{Proof of Theorem~\ref{main2} for the case $\chi_{A}\ne 0$}

It is important to note that Lenzing--Meltzer (\cite[Proposition~4.2]{LM}) shows that, if $\chi_{A}\ne 0$, 
\begin{equation}\label{standard-auto}
{\rm Auteq}(\D^b(\PP^{1}_{A,\Lambda}))\simeq({\rm Aut}(\PP^{1}_{A,\Lambda})\ltimes {\rm Pic}(\PP^{1}_{A,\Lambda}))\times\ZZ[1].
\end{equation}

\subsubsection{Case $\chi_{A}>0$}

Geigle--Lenzing (cf. \cite[Section~5.4.1]{GL}) gives an equivalence of triangulated categories
\begin{equation}
\D^{b}(\PP^{1}_{A,\Lambda})\cong \D^{b}(\CC\vec{\Delta}_{A}),
\end{equation} 
where $\vec{\Delta}_{A}$ is the extended Dynkin quiver below. 
\begin{table}[htb]
{\renewcommand\arraystretch{1.1}
\begin{tabular}{|c||c|c|c|c|c|} \hline
A& $(1,a_{2},a_{3})$ & $(2,2,a_{3})$ & (2,3,3) & (2,3,4) & (2,3,5) \\
\hline
$\vec{\Delta}_{A}$ & $\widetilde{A}_{a_{1},a_{2}}$&$\widetilde{D}_{a_{3}}$ & $\widetilde{E}_{6}$ & $\widetilde{E}_{7}$ & $\widetilde{E}_{8}$ \\
\hline
\end{tabular}
}
\end{table}

This equivalence with Corollary~\ref{GY-Dynkin} yields $h(F)=\log\rho(\N(F))$. 
Then, \cite[Theorem~4.2, Theorem~4.5]{MY} show that $\rho(\N(F))=1$. We have finished the proof.
\qed

\subsubsection{Case $\chi_{A}<0$}
We shall prove that $h(F)=0$ for each $F\in ({\rm Aut}(\PP^{1}_{A,\Lambda})\ltimes {\rm Pic}(\PP^{1}_{A,\Lambda}))\times\ZZ[1]$
if $\chi_A\leq0$.

Choose $\{[\O_{\PP^{1}_{A,\Lambda}}], [S_{1,1}],\dots, [S_{i,j}],\dots, [S_{r,a_{r}-1}],[S]\}$ as a basis of $\N(\PP^{1}_{A,\Lambda})$.
\begin{lem}\label{sr-auto1}
For $f\in {\rm Aut}(\PP^{1}_{A,\Lambda})$, the automorphism $\N(f^{*})$ is a composition of permutations
exchanging $[S_{i,j}]$ and $[S_{i',j}]$ for $j=1,\dots a_{i-1}$ if $a_{i}=a_{i'}$ and fixing $[\O_{\PP^{1}_{A,\Lambda}}]$
and $[S]$. In particular, we have $\rho(\N(f^{*}))=1$.
\end{lem}
\begin{pf}
This is a direct consequence of \cite[Proposition~3.1]{LM}.
Note also that $r\ge 3$ since $\chi_A\le 0$ and hence ${\rm Aut}(\PP^{1}_{A,\Lambda})$ is a finite group.
\qed 
\end{pf}

\begin{lem}\label{sr-auto2}
For $\L\in {\rm Pic}(\PP^{1}_{A,\Lambda})$, we have $\rho(\N(-\otimes \L))=1$.
\end{lem}
\begin{pf}
We may assume that $\L=\O_{\PP^{1}_{A, \Lambda}}(\vec{x}_{i})$.
By \cite[(2.5.3) and (2.5.4)]{GL}, for $i=1,\dots, r$ and $j=1,\dots, a_{i}-1$,
\begin{equation}
S\otimes \O_{\PP^{1}_{A, \Lambda}}(\vec{x})\cong S, \quad S_{i,j}\otimes \O_{\PP^{1}_{A, \Lambda}}(\vec{x})\cong S_{i,j+p_{i}}
\quad \text{for} \ \ \vec{x}=l\vec{c}+\sum_{i=1}^{r}p_{i}\vec{x}_{i}.
\end{equation}
It follows from the above isomorphisms that the representation matrix of $\N(-\otimes \L)$ in the basis
becomes an upper triangular matrix whose diagonal entries are all $1$. 
Hence, its spectral radius is equal to $1$. 
\qed 
\end{pf}
Each auto-equivalence $F$ is represented as $F=f^*(-\otimes\L)[m]$ (cf. (\ref{standard-auto})). 
Since Lemma~\ref{entropy} (v) gives $h(F)=h(f^*(-\otimes\L))$, we may assume $F=f^*(-\otimes\L)$. 
\begin{prop}\label{standard-curve}
We have 
\begin{equation}
h(f^*(-\otimes\L))=\log\rho(\N(f^*(-\otimes\L)))=0. 
\end{equation}
\end{prop}
\begin{pf}
Take $G$ and $G^{*}$ as in Proposition~\ref{generators-orb-proj-line}.
By Proposition~\ref{smooth-proper},  
\begin{equation*}
h(F)=\lim_{n\rightarrow\infty}\frac{1}{n}\log\delta'_{\D^b(\PP^{1}_{A,\Lambda})}(G,F^nG^*).
\end{equation*}
By straightforward calculation,
\begin{equation*}
F^{n}G^{*}=\L_{1}\otimes \L_{2}\otimes \cdots \otimes \L_{n} \otimes (f^{*})^{n}G^{*}, \quad \L_{k}:=(f^{*})^{k}\L.
\end{equation*}
Note that $f^{*}(G^{*})=G^{*}$ and $\deg(f^{*}\L)=\deg(\L)$ by Lemma~\ref{sr-auto1}.

Suppose that ${\rm deg}(\L)>0$. 
For $n\gg 0$ and $\vec{z}\in L_{A}$, we have
\begin{equation*}
\deg (\L_{1}\otimes \L_{2}\otimes \cdots \otimes \L_{n} \otimes \O_{\PP^{1}_{A,\Lambda}}(\vec{z}))=
n\deg(\L)+\deg(\O_{\PP^{1}_{A, \Lambda}}(\vec{z}))\gg 0.
\end{equation*}
Therefore, Proposition~\ref{hom-str-of-line} {\rm (i)} and {\rm (ii)} yield
\begin{eqnarray*} 
&&{\rm Ext}^{1}(G, \L_{1}\otimes \L_{2}\otimes \cdots \otimes \L_{n} \otimes G^{*})\\
&\cong &{\rm Hom}(\L_{1}\otimes \L_{2}\otimes \cdots \otimes \L_{n} \otimes G^{*}, G\otimes\O_{\PP^{1}_{A, \Lambda}}(\vec{\omega}))^{*}=0.
\end{eqnarray*}

Suppose that $\deg(\L)\le 0$.
We choose $\vec{z}\in L_{A}$ so that $\deg(\O_{\PP^{1}_{A,\Lambda}}(\vec{z}))\gg 0$.
The elements $G':=G\otimes \O_{\PP^{1}_{A,\Lambda}}(\vec{z})$ and 
$G'':=G^{*}\otimes \O_{\PP^{1}_{A,\Lambda}}(-\vec{z})$
are also split-generators. 
Therefore, Proposition~\ref{hom-str-of-line} {\rm (i)} yields
\begin{equation*}
{\rm Hom}(G', \L_{1}\otimes \L_{2}\otimes \cdots \otimes \L_{n} \otimes G'')=0.
\end{equation*}

Hence it follows from Proposition~\ref{hom-str-of-line} (iii), Lemma~\ref{sr-auto1} and Lemma~\ref{sr-auto2} that
\begin{eqnarray*}
h(F)&=& \lim_{n\rightarrow\infty}\frac{1}{n}\log\left|\chi(G,F^{n}G^{*})\right|\\
&\le&\log \rho(\N(F))=\log \rho(\N(f^{*}(-\otimes \L))=0.
\end{eqnarray*}
\qed
\end{pf}
To summarize, we have finished the proof of Theorem~\ref{main2} for the case $\chi_{A}<0$

\subsection{Proof of Theorem~\ref{main2} for the case $\chi_{A}=0$}

Define a homomorphism $\nu:\N(\PP^{1}_{A,\Lambda})\rightarrow \ZZ^{2}$ by $\nu([M]):=({\rm rk}([M]),{\rm \deg}([M]))$ 
and a skew symmetric bilinear form on $\ZZ^2$ by
\begin{equation}
\chi': \ZZ^{2}\times \ZZ^{2}\rightarrow \ZZ, \quad ((r_{1},d_{1}),(r_{2},d_{2}))\mapsto r_{1}d_{2}-r_{2}d_{1}.  
\end{equation}

\begin{lem}\label{RR-for-orb}
For $M_1,M_2\in \D^b(\PP^{1}_{A,\Lambda})$, we have 
\begin{equation}
\sum_{j=1}^{a}\chi(M_{1}(j\vec{\omega}),M_{2})=\chi'(\nu(M_1),\nu(M_2)).
\end{equation}
\end{lem}
\begin{pf}
It follows from \cite[Section 2.9]{GL} with $\chi_{A}=0$.
\qed
\end{pf}
Lemma~\ref{RR-for-orb} gives the following natural group homomorphism: 
\begin{equation}
\varphi: {\rm Auteq}(\D^{b}(\PP^{1}_{A,\Lambda}))\rightarrow {\rm Aut}_{\ZZ}(\ZZ^{2}, \chi')\cong{\rm SL}(2,\ZZ).
\end{equation} 
Denote by ${\rm Pic}^{0}(\PP^{1}_{A,\Lambda})\subset {\rm Pic}(\PP^{1}_{A,\Lambda})$ the subgroup consisting
of elements with degree zero.
\begin{prop}
There exists the following exact sequence:
\begin{equation}
\{1\}\rightarrow ({\rm Aut}(\PP^{1}_{A,\Lambda})\ltimes {\rm Pic}^{0}(\PP^{1}_{A,\Lambda}))\times \ZZ[2]
\rightarrow {\rm Auteq}(\D^{b}(\PP^{1}_{A,\Lambda})) \xrightarrow{\varphi} {\rm SL}(2,\ZZ)\rightarrow \{1\}.
\end{equation}
\end{prop}
\begin{pf}
This is a direct consequence of \cite[Theorem~6.3]{LM}.
\qed
\end{pf}

\begin{lem}\label{ker-phi}
The map $h:{\rm Auteq}(\D^b(\PP^{1}_{A,\Lambda}))\longrightarrow \RR_{\ge 0}$, $F\mapsto h(F)$ factors 
through ${\rm SL}(2,\ZZ)$.
\end{lem}
\begin{pf}
Choose $\O_{\PP^{1}_{A,\Lambda}}(\vec{z})$ so that $\deg(\O_{\PP^{1}_{A,\Lambda}}(\vec{z}))\gg 0$
and set $G':=G\otimes \O_{\PP^{1}_{A,\Lambda}}(\vec{z})$, $G'':=G^{*}\otimes \O_{\PP^{1}_{A,\Lambda}}(-\vec{z})$.
By Lemma~\ref{entropy} (v), we can assume that 
an element $F\in {\rm Auteq}(\D^b(\PP^{1}_{A,\Lambda}))$ is of the form $F=F'F_1$ with $F'\in {\rm Auteq}(\D^b(\PP^{1}_{A,\Lambda}))$ and $F_1\in {\rm Aut}(\PP^{1}_{A,\Lambda})\ltimes {\rm Pic}^{0}(\PP^{1}_{A,\Lambda})$. 
Then there exist $F_2,\dots, F_n\in {\rm Aut}(\PP^{1}_{A,\Lambda})\ltimes {\rm Pic}^{0}(\PP^{1}_{A,\Lambda})$ such that 
$F^n=(F'F_1)^n=F'^nF_n \cdots  F_1$. We have 
\[
\delta_{\D^b(\PP^{1}_{A,\Lambda})}(G',F^{n} G'')
\le \delta_{\D^b(\PP^{1}_{A,\Lambda})}(G',F'^n G')\delta_{\D^b(\PP^{1}_{A,\Lambda})}(G',F_n\cdots F_1 G''),
\]
and hence, 
\begin{equation*}
h(F)\le h(F')+\lim_{n\rightarrow\infty}\frac{1}{n}\log\delta_{\D^b(\PP^{1}_{A,\Lambda})}(G',F_n\cdots F_1 G'').
\end{equation*}
The functor $F_i$ is of the form $f_i^*(-\otimes\L_i)$ for some $f_i\in {\rm Aut}(\PP^{1}_{A,\Lambda})$ and
$\L_i\in {\rm Pic}^{0}(\PP^{1}_{A,\Lambda})$. 
For arbitrary $\vec{y}\in L_{A}$, we have 
$\deg(F_n\cdots F_1(\O_{\PP^{1}_{A,\Lambda}}(\vec{y})))=\deg(\O_{\PP^{1}_{A,\Lambda}}(\vec{y}))$.
Therefore, it follows from Proposition~\ref{hom-str-of-line} (i),(iii) that 
$\delta'_{\D^b(\PP^{1}_{A,\Lambda})}(G',F_n\cdots F_1 G'')=\left|\chi(G',F_n\cdots F_1 G'')\right|$. 
Lemma~\ref{delta-delta'}, Lemma~\ref{sr-auto1} and Lemma~\ref{sr-auto2} yield
\begin{eqnarray*}
\lim_{n\rightarrow\infty}\frac{1}{n}\log\delta_{\D^b(\PP^{1}_{A,\Lambda})}(G', F_n\cdots F_1 G'')
&=&\lim_{n\rightarrow\infty}\frac{1}{n}\log\delta'_{\D^b(\PP^{1}_{A,\Lambda})}(G', F_n\cdots F_1 G'')\\
&=&\lim_{n\rightarrow\infty}\frac{1}{n}\log\left|\chi(G',F_n\cdots F_1G'')\right|=0,
\end{eqnarray*}
and hence $h(F)\leq h(F')$. 
We also have $h(F')\leq h(F)$ since $F'=FF_1^{-1}$ and $F_1^{-1}$ belongs to  ${\rm Aut}(\PP^{1}_{A,\Lambda})\ltimes {\rm Pic}^{0}(\PP^{1}_{A,\Lambda})$. 
\qed
\end{pf}

\begin{prop}\label{GY-orb}
We have $h(F)=\log\rho(\N(F))$. 
\end{prop}
\begin{pf}
Since $h(F[1])=h(F)$, we may assume that ${\rm tr}(\varphi(F))\ge 0$.
It is easy to calculate $h(F)=0$ if ${\rm tr}(\varphi(F))=0,1$ since $\varphi(F)$ is of finite order and hence $F$ is of finite order up to ${\rm Aut}(\PP^{1}_{A,\Lambda})\ltimes {\rm Pic}^{0}(\PP^{1}_{A,\Lambda})$.
If ${\rm tr}(\varphi(F))=2$, then $F=(-\otimes \O_{\PP^{1}_{A,\Lambda}}(\vec{x}))F'$ with $F'\in {\rm Aut}(\PP^{1}_{A,\Lambda})\ltimes {\rm Pic}^{0}(\PP^{1}_{A,\Lambda})$ for some $\vec{x}\in L_A$.
It follows from Proposition~\ref{standard-curve} that 
\begin{eqnarray*}
h(F)=h(-\otimes \O_{\PP^{1}_{A,\Lambda}}(\vec{x}))=0=\log\rho(\N(F)).
\end{eqnarray*}

Suppose now that ${\rm tr}(\varphi(F))>2$.
\begin{lem}\label{slope-criterion}
For indecomposable objects $M_{1}, M_{2}\in {\rm Coh}(\PP^{1}_{A,\Lambda})$, 
we have 
\begin{equation}
{\rm Ext}^1(M_{1},M_{2})=0\quad \text{ if }\quad \chi'(\nu(M_1),\nu(M_2))>0.
\end{equation}
\end{lem}
\begin{pf}
The statement follows from the slope-stability for orbifold projective lines (Proposition~5.2 in \cite{GL}) and Proposition~\ref{hom-str-of-line}~(ii).
\qed
\end{pf}

\begin{lem}[Proposition~4.6 in \cite{Kik}]\label{Z-conjugate}
Assume that ${\rm tr}(\varphi(F))>2$. 
There exists a sequence of positive integers ${\bf m}=(m_{2n},\dots,m_1)$ with $n\ge 1$
such that $\varphi(F)$ is conjugate in ${\rm SL}(2,\ZZ)$ to
\begin{equation}
\begin{pmatrix}
1& m_{2n-1}\\
m_{2n}& 1+m_{2n-1}m_{2n}
\end{pmatrix}
\cdots 
\begin{pmatrix}
1& m_{1}\\
m_{2}& 1+m_{1}m_{2}
\end{pmatrix}.
\end{equation}
\end{lem}

For each sequence of positive integers ${\bf m}=(m_{2n},\dots,m_1)$ with $n\ge 1$, 
take $F_{\bf m}\in {\rm Auteq}(\D^b(\PP^{1}_{A,\Lambda}))$ so that 
\[
\varphi(F_{\bf m})=\begin{pmatrix}
1& m_{2n-1}\\
m_{2n}& 1+m_{2n-1}m_{2n}
\end{pmatrix}
\cdots 
\begin{pmatrix}
1& m_{1}\\
m_{2}& 1+m_{1}m_{2}
\end{pmatrix}.
\]
For positive $\vec{x},\vec{y}\in L_{A}$, an elementary calculation gives
\begin{equation*}
\chi'(\nu(\O_{\PP^{1}_{A,\Lambda}}(-\vec{x})), \nu(F_{\bf m}^{n}\O_{\PP^{1}_{A,\Lambda}}(\vec{y})))>0.
\end{equation*}

Take $G, G^{*}$ as in Proposition~\ref{generators-orb-proj-line} and a positive $\vec{x}\in L_A$.
By Proposition~\ref{hom-str-of-line} {\rm (iii)}, Lemma \ref{slope-criterion}, we obtain
\begin{eqnarray*}
h(F_{\bf m})&=&\lim_{n\rightarrow\infty}\frac{1}{n}\log\delta'_{\D^b(\PP^{1}_{A,\Lambda})}(G^*\otimes\O_{\PP^{1}_{A, \Lambda}}(-\vec{x}),F_{\bf m}^{n}(G\otimes\O_{\PP^{1}_{A, \Lambda}}(\vec{x})))\\
&=&\lim_{n\rightarrow\infty}\frac{1}{n}\log\left|\chi(G^*\otimes\O_{\PP^{1}_{A, \Lambda}}(-\vec{x}),F_{\bf m}^{n}(G\otimes\O_{\PP^{1}_{A, \Lambda}}(\vec{x})))\right|\leq\log\rho(\N(F_{\bf m})). 
\end{eqnarray*}

\begin{lem}\label{alg-integer}
We have
\begin{equation}
\rho(\N(F))=\rho(\varphi(F)).
\end{equation}
In particular, $\rho(\N(F))$ is an algebraic number. 
\end{lem}
\begin{pf}
The inequality $\rho(\N(F))\geq\rho(\varphi(F))$ follows from the commutativity: $\varphi(F)\circ\nu=\nu\circ\N(F)$. 
The fact that $\varphi$ factors thorough the surjection ${\rm Aut}_\ZZ(\N(\PP^{1}_{A, \Lambda}),\chi)\to{\rm Aut}_{\ZZ}(\ZZ^{2}, \chi')$ (\cite[Theorem~7.3]{LM}) yields the reversed inequality. 
\qed
\end{pf}
Since $\varphi(F_{{\bf m}})$ is conjugate to $\varphi(F)$, it follows from Lemma~\ref{ker-phi} and Lemma~\ref{alg-integer} that
\begin{equation*}
h(F)=h(F_{{\bf m}})\leq\log\rho(\N(F_{{\bf m}}))=\log\rho(\varphi(F_{{\bf m}}))=\log\rho(\varphi(F))=\log\rho(\N(F)).
\end{equation*}
By Theorem~\ref{nGro-lower}, we finished the proof of Proposition~\ref{GY-orb}, hence of Theorem~\ref{main2}. 
\qed
\end{pf}


\end{document}